\documentclass[11pt]{article}
\usepackage{mathrsfs}
\usepackage[centertags]{amsmath}
\usepackage{amsfonts}
\usepackage{amssymb}
\usepackage{amsthm}
\usepackage{cases}
\usepackage{indentfirst}
\usepackage{epsfig}
\usepackage{graphicx}
\usepackage{color}
\usepackage[title,titletoc,header]{appendix}

\date{}

\definecolor{c20}{rgb}{0.,0.7,0.}
\definecolor{c30}{rgb}{0.,0.,1.}
\definecolor{c40}{rgb}{1,0.1,0.7}
\definecolor{c50}{rgb}{1,0,0}
\definecolor{c60}{rgb}{0.00,0.00,1.00}

\def\zw#1{\textcolor{c50}{#1}}
\def\zw#1{#1}

\def\pzx#1{\textcolor{c60}{#1}}
\def\pzx#1{#1}

\newtheorem{theorem}{Theorem}[section]

\newtheorem{lemma}{Lemma}[section]

\newtheorem{remark}{Remark}[section]

\numberwithin{equation}{section}

\newcommand\abs[1]{\left\lvert #1 \right\rvert}
\def\P{\operatorname*{\mathbb{P}}}
\def\E{\operatorname*{E}}
\def\Var{\operatorname*{Var}}
\def\Cov{\operatorname*{Cov}}

\def\R{\operatorname*{\mathbb{R}}}

\allowdisplaybreaks[4]
\usepackage{setspace}
\begin{document}
\title{Asymptotic behaviors of bivariate Gaussian powered extremes }
\author{ {Wei Zhou\quad and\quad  Zuoxiang Peng\footnote{Corresponding author. Email: pzx@swu.edu.cn}}\\
{\small School of Mathematics and Statistics, Southwest University,
Chongqing, 400715, China}}
\date{}

\maketitle
\begin{quote}
{\bf Abstract}~~ In this paper, joint asymptotics of powered maxima for a triangular array of bivariate powered Gaussian random vectors are considered. Under the H\"usler-Reiss condition, limiting distributions of powered maxima are derived. Furthermore, the second-order expansions of the joint distributions of powered maxima are established under the refined
H\"usler-Reiss condition.


\end{quote}

\section{Introduction}
\label{sec1}

For independent and identically distributed bivariate Gaussian random vectors  with constant coefficient in each vector, Sibuya (1960) showed that componentwise maxima are asymptotically independent, and Embrechts et al. (2003) proved the asymptotical independence in the upper tail. To overcome those shortcomings in its applications, H\"usler and Reiss (1989) considered the asymptotic behaviors of extremes of Gaussian triangular arrays with varying coefficients. Precisely, let $\{(X_{ni},Y_{ni}), 1\leq i\leq n, n\geq 1\}$ be a
triangular array of independent bivariate Gaussian random vectors with $\E X_{ni}=\E Y_{ni}=0$, $\Var X_{ni}=\Var Y_{ni}=1$ for  $1\le i\le n$, $n\ge 1$.
and $\Cov(X_{ni},Y_{ni})=\rho_n$. Let $F_{\rho_{n}}(x,y)$ denote the joint distribution of vector $(X_{ni},Y_{ni})$ for $i\le n$. The partial maxima $\mathbf{M_n}$ is defined by
\begin{eqnarray*}
\mathbf{M_n}=(M_{n1},M_{n2})=(\max_{1\leq i\leq n}X_{ni},\max_{1\leq i\leq n}Y_{ni}).
\end{eqnarray*}
 H\"usler and Reiss (1989) and  Kabluchko (2009) showed that
 \begin{eqnarray}
\label{eq1.1}
\lim_{n\to\infty}\P\left( M_{n1}\leq b_n+\frac{x}{b_n},M_{n2}\leq b_n+\frac{y}{b_n} \right)= H_\lambda(x,y)
\end{eqnarray}
holds if and only if the following H\"usler-Reiss condition
\begin{eqnarray}
\label{eq1.2}
\lim_{n\to \infty} b_n^2(1-\rho_n) = 2\lambda^2 \in [0,\infty]
\end{eqnarray}
holds, where the normalizing constant $b_n$ satisfying
\begin{eqnarray}
\label{eq1.3}
1-\Phi(b_n)=n^{-1}
\end{eqnarray}
and the max-stable H\"usler-Reiss distribution is given by
\begin{eqnarray}
\label{eq1.4}
H_{\lambda}(x,y)=\exp\left(-\Phi\left(\lambda+\frac{x-y}{2\lambda}\right)e^{-y}-
\Phi\left(\lambda+\frac{y-x}{2\lambda}\right)e^{-x}\right), \quad x,y\in \mathbb{R},
\end{eqnarray}
where $\Phi(\cdot)$ and $\varphi(\cdot)$ denote respectively the distribution function and density function of a standard Gaussian random variable.
Note that $H_0(x,y)=\Lambda(\min(x,y))$ and $H_{\infty}(x,y)=\Lambda(x)\Lambda(y)$ with $\Lambda(x)=\exp(-\exp(-x))$, $x\in\R$.

Recently, contributions to H\"usler-Reiss distribution and its extensions
are achieved considerably. For instance, Hashorva (2005, 2006) showed that the limit distributions
of maxima also holds for triangular arrays of general bivariate elliptical distributions
if the distribution of random radius is in the Gumbel or Weibull max-domain of attraction, and Hashorva and Ling (2016) extended the results to bivariate skew elliptical triangular arrays. For more work on asymptotics of bivariate triangular arrays, see Hashorva (2008, 2013) and Hashorva et al. (2012).

Higher-order expansions of distributions of extremes on  H\"usler-Reiss bivariate Gaussian triangular arrays were considered firstly by Hashorva et al. (2016)
provided that $\rho_{n}$ satisfies the following refined H\"usler-Reiss condition
 \begin{eqnarray}
\label{eq1.5}
\lim_{n\to\infty}b_n^2 (\lambda_n-\lambda) = \alpha  \in \R,
\end{eqnarray}
where $\lambda_n=(\frac{1}{2}b_n^2(1-\rho_n))^{1/2}$ and $\lambda \in(0,\infty) $,
with $b_n$ given by \eqref{eq1.3}.
Uniform convergence rate was considered by Liao and Peng (2014).
For copula version of the limit in  H\"usler-Reiss model, Frick and Reiss (2013)
 considered the penultimate and ultimate convergence
rates for distribution of $( n(\max_{1\leq i\leq n}\Phi(X_{ni})-1),
 n((\max_{1\leq i\leq n}\Phi(Y_{ni})-1)) )$, and  Liao et al. (2016)
extended the results to the settings of $n$ independent and non-identically distributed observations,
where the $i$th observation follows from normal copula with correlation coefficient
 being either a parametric or a nonparametric function of $i/n$.

The objective of this paper is to study the asymptotics of powered-extremes of H\"usler-Reiss bivariate Gaussian triangular arrays. Interesting results in Hall (1980) showed that the convergence rates of the distributions of powered-extremes of  independent and identically distributed univariate Gaussian sequence depend on the power index and normalizing constants. Precisely, Let $\abs{M_n}^t$ denote the powered maximum with any power index $t>0$, then
 \begin{eqnarray}
\label{eq1.6}
\lim_{n \to \infty} b_n^{2} \Big[\P\left(\abs{M_n}^t\leq
c_n x+{d_n}\right)-\Lambda(x)\Big]
=\Lambda(x)\pzx{\mu}(x)
\end{eqnarray}
with normalizing constants $c_{n}$ and $d_{n}$ given by
\begin{eqnarray}
\label{eq1.7}
c_n=t b_n^{t-2}, \quad d_n=b_n^t, \quad t>0.
\end{eqnarray}
Furthermore, for $t=2$ with normalizing constants $c_{n}^{*}$ and $d_{n}^{*}$ given by
\begin{eqnarray}
\label{add2}
c_{n}^{*}=2 - 2b_n^{-2}, \quad d_{n}^{*}=b_{n}^{2} - 2b_n^{-2},
\end{eqnarray}
we have
 \begin{eqnarray}
\label{add1}
\lim_{n \to \infty} b_n^{4} \Big[\P\left(\abs{M_n}^2\leq
c_{n}^{*} x+d_{n}^{*}\right)-\Lambda(x)\Big]
=\Lambda(x) \pzx{\nu}(x),
\end{eqnarray}
where $b_n$ is defined in \eqref{eq1.3}, and $\mu(x)$ and $\nu(x)$ are respectively given by
\begin{eqnarray}
\label{eq1.8}
\mu(x)=\left(1+x+\frac{2-t}{2}x^2\right) e^{-x}, \quad
\nu(x)=-\left( \frac{7}{2}+3x+x^2\right)e^{-x}.
\end{eqnarray}

Motivated by findings of H\"usler-Reiss (1989), Hall (1980) and Hashorva et al. (2016), we will consider the distributional asymptotics of  powered-extremes \zw{of} H\"usler-Reiss bivariate Gaussian triangular arrays, and hope that the convergence rates can be improved as $t=2$, similar to \eqref{add1} in univariate case. Unfortunately, our results provide negative answers except
\zw{two extreme cases}.

The rest of the paper is organized as follows. In Section \ref{sec2} we provide
the main results and all proofs are deferred to Section \ref{sec4}. Some auxiliary
results are given in Section \ref{sec3}.

\section{Main Results}
\label{sec2}

In this section, the limiting distributions and the second-order expansions on distributions of
normalized bivariate powered-extremes are provided if $\rho_n$ satisfies \eqref{eq1.2} and \eqref{eq1.5}, respectively.
The first main result, stated as follows, is the limit distributions of bivariate normalized powered-extremes.

\begin{theorem}\label{thm1}
Let the norming constants $c_n$ and $d_n$ be given by \eqref{eq1.7}. Assume that \eqref{eq1.2} holds with
$\lambda \in (0,\infty)$. Then for all $x,y \in \R$, we have
\begin{eqnarray}
\label{eq2.2}
\lim_{n\to \infty}\P\left( \abs{M_{n1}}^t\leq c_n x+d_n, \abs{M_{n2}}^t\leq c_n y+d_n \right)
=H_{\lambda}(x,y).
\end{eqnarray}
\end{theorem}

\begin{remark}\label{remark1}
For $t=2$,  with arguments similar to the proof of Theorem \ref{thm1} one can show that \eqref{eq2.2} also holds with $c_{n}$ and $d_{n}$ being replaced by $c_{n}^{*}$ and $d_{n}^{*}$ given by \eqref{add2}.
\end{remark}

Next we investigate the convergence rate of
\begin{equation}\label{add3}
\Delta(F_{\rho_n,t}^{n},H_{\lambda};\pzx{c_{n}, d_{n}}; x,y)=
\P\left( \abs{M_{n1}}^t\leq c_n x+d_n, \abs{M_{n2}}^t\leq c_n y+d_n \right)
-H_{\lambda}(x,y)\to 0
\end{equation}
as $n\to\infty$ under the refined second-order  H\"usler-Reiss condition \eqref{eq1.5}. The results are stated as follows.

\begin{theorem}\label{thm2}
If the second H\"usler-Reiss condition \eqref{eq1.5} holds with
$\lambda_n=(\frac{1}{2}b_n^2(1-\rho_n))^{1/2}$ and $\lambda \in (0,\infty)$, then for all $x,y \in \R$, we have
\begin{eqnarray}
\label{eq2.3}
\lim_{n\to \infty}(\log n)\Delta(F_{\rho_n,t}^{n},H_{\lambda};\pzx{c_{n}, d_{n}};x,y)
=\frac{1}{2} \tau(\alpha,\lambda,x,y,t) H_{\lambda}(x,y)
\end{eqnarray}
with
\begin{eqnarray*}
  \tau(\alpha,\lambda,x,y,t)&=& \mu(x)\Phi\left(\lambda+\frac{y-x}{2\lambda}\right)+\mu(y)\Phi\left(\lambda+\frac{x-y}{2\lambda}\right)
 \\
 && +\Big(2\alpha-(x+y+2)\lambda-\lambda^3\Big)e^{-x}\varphi\left(\lambda+\frac{y-x}{2\lambda}\right),
\end{eqnarray*}
where $\mu(x)$ is the one given by \eqref{eq1.8}.
\end{theorem}

\begin{remark}\label{remark2}
For $t=2$,  let $c_n$ and $d_n$ be replaced by $c_n^\ast$ and
 $d_n^\ast$ respectively in \eqref{add2}, one can show that
 \begin{eqnarray*}
\lim_{n\to \infty}(\log n)\Delta(F_{\rho_n,2}^{n},H_{\lambda};\pzx{c_{n}^{*}, d_{n}^{*}};x,y)
=\frac{1}{2} \chi(\alpha,\lambda,x,y) H_{\lambda}(x,y)
\end{eqnarray*}
with
\begin{eqnarray*}
\chi(\alpha,\lambda,x,y)
= \Big(2\alpha-(x+y+2)\lambda-\lambda^3\Big)e^{-x}\varphi\left(\lambda+\frac{y-x}{2\lambda}\right).
\end{eqnarray*}
The result shows the fact that the convergence rates can not be improved as $t=2$ with normalizing constants $c_{n}^{*}$ and $d_{n}^{*}$, contrary to the result of univariate Gaussian case provided by Hall (1980).
\end{remark}

In order to obtain the convergence rates of \eqref{add3} for two extreme cases $\lambda=0$ and $\lambda=\infty$, we may need some additional conditions.
Following results show that rates of convergence are considerably different with different choice of \pzx{normalizing} constants. With power index $t>0$ and \pzx{normalizing} constants $c_{n}$ and $d_{n}$ given by \eqref{eq1.7}, the results are stated as follows.

\begin{theorem}\label{thm3}
Let $c_n$ and $d_n$ be given by \eqref{eq1.7}. With  $x,y \in \R$ and $t>0$ we have the following results.
\begin{itemize}
\item[(a).]
For the case of $\lambda=\infty$,

(i) if $\rho_n \in [-1,0]$, we have
\begin{eqnarray}\label{eq2.4}
\lim_{n\to \infty}(\log n)\Delta(F_{\rho_n,t}^{n},H_{\infty};\pzx{c_{n}, d_{n}};x,y)
=\frac{1}{2}\Big(\mu(x)+\mu(y)\Big)  H_{\infty}(x,y).
\end{eqnarray}

(ii) if $\rho_n \in (0,1)$ and $\lim_{n \to \infty} \frac{\log b_n}{b_n^2(1-\rho_n)}=0$,
then \eqref{eq2.4} also holds.

\item[(b).]
For the case of $\lambda=0$,

(i) if $\rho_n=1$, we have
\begin{eqnarray}\label{eq2.5}
\lim_{n\to \infty}(\log n)\Delta(F_{\rho_n,t}^{n},H_0;\pzx{c_{n}, d_{n}};x,y)
=\frac{1}{2} \mu\Big(\min(x,y)\Big)  H_0(x,y).
\end{eqnarray}
(ii) if $\rho_n \in (0,1)$ and $\lim_{n \to \infty}b_n^{6}(1-\rho_n)=c_{1}\in [0,\infty)$, then
\eqref{eq2.5} also holds.
\end{itemize}
\end{theorem}

Theorem \ref{thm3} shows that convergence rates of \eqref{add3} are the same order of $1/\log n$ if we choose the normalizing constants $c_{n}$ and $d_{n}$ given by \eqref{eq1.7} as $t>0$. With another pair of normalizing constants $c_{n}^{*}$ and $d_{n}^{*}$ given by \eqref{add2}, the following results show that convergence rates of \eqref{add3} can be improved.
\begin{theorem}\label{thm4}
For $t=2$, let $c_{n}^{*}$ and $d_{n}^{*}$ be given by \eqref{add2}. With $x,y \in \R$ we have the following results.
\begin{itemize}
\item[(a).]
For the case of $\lambda=\infty$,

(i) if $\rho_n \in [-1,0]$, we have
\begin{eqnarray}\label{eq2.6}
\lim_{n\to \infty}(\log n)^2\Delta(F_{\rho_n,2}^{n},H_{\infty};c_{n}^{*}, d_{n}^{*};x,y)
=\frac{1}{4} \Big(\nu(x)+\nu(y)\Big) H_{\infty}(x,y).
\end{eqnarray}

(ii) if $\rho_n \in (0,1)$ and $\lim_{n \to \infty} \frac{\log b_n}{b_n^2(1-\rho_n)}=0$,
then \eqref{eq2.6} also holds.

\item[(b).]
For the case of $\lambda=0$,

(i) if $\rho_n=1$, we have
 \begin{eqnarray}\label{eq2.7}
\lim_{n\to \infty}(\log n)^2\Delta(F_{\rho_n,2}^{n},H_0;c_{n}^{*}, d_{n}^{*};x,y)
=\frac{1}{4}\nu\Big(\min(x,y)\Big) H_0(x,y).
\end{eqnarray}

(ii) if $\rho_n \in (0,1)$ and $\lim_{n \to \infty}b_n^{14}(1-\rho_n)=c_2\in [0,\infty)$, then
\eqref{eq2.7} also holds.
\end{itemize}

\end{theorem}

\section{Auxiliary Lemmas}
\label{sec3}
For notational simplicity, let
\begin{eqnarray}
\label{eq2.1}
\omega_{n,t}(x)=(c_n x+d_n)^{1/t}
\quad
\mbox{for} \quad t>0,\quad \mbox{and}\quad \omega_{n,2}^{*}(x)=(c_{n}^{*} x+d_{n}^{*})^{1/2}
\quad
\mbox{as} \quad t=2,
\end{eqnarray}
where the normalizing constants $c_n$ and $d_n$, and $c_{n}^{*}$ and $d_{n}^{*}$ are those given by \eqref{eq1.7} and \eqref{add2}, respectively. Define
\begin{eqnarray*}
\pzx{\bar{\Phi}(z)=1-\Phi(z)},\quad\quad
\zw{\bar{\Phi}_{n.t}(z)= n\bar{\Phi}(\omega_{n,t}(z))}
\end{eqnarray*}
and
\begin{eqnarray}
\label{eq3.1a}
I_k:=\int_{y}^{\infty}\varphi\left(\lambda+\frac{x-z}{2\lambda}\right)e^{-z}z^k dz,
\quad k=0,1,2.
\end{eqnarray}

\begin{lemma}\label{lemma1}
Under  the conditions of Theorem \ref{thm1}, we have
\begin{eqnarray}
\label{eq3.1}
\lim_{n\to \infty} \P\Big( M_{n1}\leq \omega_{n,t}(x),M_{n2}\leq \omega_{n,t}(y) \Big)
=H_{\lambda}(x,y)
\end{eqnarray}
\end{lemma}

\noindent
\textbf{Proof.}~~ With the choice of $c_n$ and $d_n$ in \eqref{eq1.7}, it follows from \eqref{eq2.1}
that
\begin{eqnarray*}
\omega_{n,t}(z)=(c_n z+d_n)^{1/t}=b_n \left( 1+ \zw{ z b_n^{-2}+ \frac{1-t}{2}z^2 b_n^{-4}} +O(b_n^{-6}) \right)
\end{eqnarray*}
for fixed $z$, hence for fixed $x$ and $z$,
\begin{eqnarray}\label{eq3.1b}
\nonumber
 &&
\frac{\omega_{n,t}(x)-\rho_n \omega_{n,t}(z)}{\sqrt{1-\rho_n^2}}
\\ \nonumber
&=& b_n \sqrt{\frac{1-\rho_n}{1+\rho_n}}+\frac{x-z}{b_n\sqrt{1-\rho_n^2}}
+\frac{z}{b_n}\sqrt{\frac{1-\rho_n}{1+\rho_n}}
+\frac{(1-t)(x^2-z^2)}{2b_n^3 \sqrt{1-\rho_n^2}}
+\frac{(1-t)z^2}{2b_n^3}\sqrt{\frac{1-\rho_n}{1+\rho_n}}
\\
\nonumber
&&
+\sqrt{\frac{1-\rho_n}{1+\rho_n}}O(b_n^{-5})
\\
&=& \left( \lambda_n+\frac{x-z}{2\lambda_n} \left( 1+\frac{(1-t)(x+z)}{2 b_n^2} \right) +\frac{\lambda_n z}{b_n^2} +
\frac{(1-t)\lambda_n z^2}{2b_n^{4}} +
\lambda_n O(b_n^{-6})\right)(1-\frac{\lambda_n^2}{b_n^2})^{-\frac{1}{2}}
\end{eqnarray}
which implies that
\begin{eqnarray}\label{add4}
\lim_{n \to \infty}\frac{\omega_{n,t}(x)-\rho_n \omega_{n,t}(z)}{\sqrt{1-\rho_n^2}}=\lambda+\frac{x-z}{2\lambda}
\end{eqnarray}
holds since $\lambda_n\to \lambda$ as $n \to \infty$.

With $a_{n}=1/b_{n}$ it follows from \eqref{add4} that
\begin{eqnarray}
\label{eq3.2}
&&
n \P \left( X>\omega_{n,t}(x),Y>\omega_{n,t}(y) \right)
\nonumber\\
&=& n \int_{\omega_{n,t}(y)}^{\infty}\bar{\Phi}\left(\frac{\omega_{n,t}(x)-\rho_n z}{\sqrt{1-\rho_n^2}}\right)d\Phi(z)
\nonumber\\
&=& \frac{b_n}{\varphi(b_n)} (1-b_n^{-2}+O(b_n^{-4}))^{-1} \int_{y}^{\infty}\bar{\Phi}\left(\frac{\omega_{n,t}(x)-\rho_n \omega_{n,t}(z)}{\sqrt{1-\rho_n^2}}\right)\varphi(b_n(1+\pzx{tz}a_n^2)^{1/t})d(b_n(1+\pzx{tz}a_n^2)^{1/t})
\nonumber\\
&=&  (1+b_n^{-2}+O(b_n^{-4})) \int_{y}^{\infty}\bar{\Phi}\left(\frac{\omega_{n,t}(x)-\rho_n \omega_{n,t}(z)}{\sqrt{1-\rho_n^2}}\right)
\exp\left(\frac{b_n^2}{2}\left( 1-(1+tza_n^2)^{2/t}\right)\right) (1+tza_n^2)^{1/t-1}dz
\nonumber\\
&\to& \int_{y}^{\infty}\bar\Phi\left(\lambda+\frac{x-z}{2\lambda}\right)e^{-z}dz
\nonumber\\
&=& e^{-y}+e^{-x}-\Phi\left(\lambda+\frac{x-y}{2\lambda}\right)e^{-y}-
\Phi\left(\lambda+\frac{y-x}{2\lambda}\right)e^{-x}
\end{eqnarray}
as $n\to\infty$. Meanwhile, one can check that
\begin{eqnarray}\label{add8}
\lim_{n \to \infty}\zw{\bar{\Phi}_{n,t}(x)}=e^{-x}.
\end{eqnarray}
It follows from \eqref{eq3.2} and \eqref{add8} that
\begin{eqnarray*}
&&   \P\left( M_{n1}\leq \omega_{n,t}(x),M_{n2}\leq \omega_{n,t}(y) \right)
\\
&=& \exp \Big[ -\zw{\bar{\Phi}_{n,t}(x)-\bar{\Phi}_{n,t}(y)}
 + n \P \left( X>\omega_{n,t}(x),Y>\omega_{n,t}(y) \right) +o(1) \Big]
\\
&\to & H_{\lambda}(x,y)
\end{eqnarray*}
as $n\to\infty$. The desired result follows. \qed

Following result is useful to the proof of Lemma \ref{lemma2}.
\begin{lemma}\label{app}
With $a_{n}=1/b_{n}$, for large $n$ we have
\begin{eqnarray}\label{add5}
\nonumber
&& \int_{y}^{\infty} \Phi\left(\frac{\omega_{n,t}(x)-\rho_n \omega_{n,t}(z)}{\sqrt{1-\rho_n^2}}\right)
 \exp\left(\frac{b_n^2}{2}\left( 1-(1+tza_n^2)^{2/t}\right)\right) (1+tza_n^2)^{1/t-1}\,dz\\
&=&
\int_{y}^{\infty}   \Phi\left(\frac{\omega_{n,t}(x)-\rho_n \omega_{n,t}(z)}{\sqrt{1-\rho_n^2}}\right)
\left( 1+\left((1-t)z-\frac{2-t}{2}z^2\right)\zw{b_n^{-2}} \right)e^{-z} dz+O(b_n^{-4}).
\end{eqnarray}
\end{lemma}

\noindent
\textbf{Proof.}~~First note for large $n$ and $\abs{x}\leq \frac{b_n^2}{4(4+t)}$,
\begin{eqnarray}
\nonumber
\label{bounded}
&& \abs{ \exp\left(\frac{b_n^2}{2}\left( 1-(1+\pzx{tx} a_n^2)^{2/t}\right)\right) (1+\pzx{tx} a_n^2)^{1/t-1}
-e^{-x}\left( 1+\frac{1}{b_n^2}\left((1-t)x-\frac{2-t}{2}x^2\right) \right)}
\\
&\leq& b_n^{-4}s(x)\exp \left(-x+\frac{\abs{x}}{4} \right),
\end{eqnarray}
where $a_{n}=1/b_{n}$ and $s(x)\ge 0$ ia a polynomial on $x$ independent of $n$,
cf.\zw{ Lemma 3.2} in Li and Peng (2016).

It follows from \eqref{bounded} that
\begin{eqnarray}
\nonumber
\label{more1}
&& \int_{y}^{4 \log b_n} \abs{ \Phi\left(\frac{\omega_{n,t}(x)-\rho_n \omega_{n,t}(z)}{\sqrt{1-\rho_n^2}}\right)
 \exp\left(\frac{b_n^2}{2}\left( 1-(1+\pzx{tz} a_n^2)^{2/t}\right)\right) (1+\pzx{tz} a_n^2)^{1/t-1} \right.
 \\
 \nonumber
 && \quad \left.- \Phi\left(\frac{\omega_{n,t}(x)-\rho_n \omega_{n,t}(z)}{\sqrt{1-\rho_n^2}}\right)
e^{-z}\left( 1+b_n^{-2}\left((1-t)z-\frac{2-t}{2}z^2\right) \right)    } dz
\\
\nonumber
&\leq & \int_{y}^{4 \log b_n} \abs{\Phi\left(\frac{\omega_{n,t}(x)-\rho_n \omega_{n,t}(z)}{\sqrt{1-\rho_n^2}}\right)}
b_n^{-4}s(z)\exp \left(-\frac{3z}{4} \right)dz
\\
\nonumber
&<& b_n^{-4} \int_{y}^{4 \log b_n}s(z)\exp \left(-\frac{3z}{4}\right)dz
\\
&=& O(b_n^{-4})
\end{eqnarray}
and
\begin{eqnarray}
\label{more5}
\nonumber
&& \int_{ 4 \log b_n}^{\infty} e^{-\frac{z}{2}} \left( e^{-\frac{z}{2}}\left(
1+b_n^{-2}\left(\abs{1-t}z+\frac{\abs{2-t}}{2}z^2\right) \right) \right) dz
\\
\nonumber
&\leq & e^{-2\log b_n} \left (
1+\zw{b_n^{-2}}\left(4\abs{1-t}\log b_n+8\abs{2-t}(\log b_n)^2\right) \right)\int_{ 4 \log b_n}^{\infty}e^{-\frac{z}{2}}dz
\\
\nonumber
&=& 2b_n^{-4} \left(
1+b_n^{-2}\left(4\abs{1-t}\log b_n+8\abs{2-t}(\log b_n)^2\right) \right)
\\
&=& O(b_n^{-4}).
\end{eqnarray}
So, the remainder is to show
\begin{eqnarray}
\label{more2}
 A_{n}=\int_{4 \log b_n}^{\infty}
 \exp\left(\frac{b_n^2}{2}\left( 1-(1+\pzx{tz} a_n^2)^{2/t}\right)\right) (1+\pzx{tz} a_n^2)^{1/t-1}dz
= O(b_n^{-4})
\end{eqnarray}
for large $n$. We check \eqref{more2} in turn for
\zw{ $0<t< 1$ and $t\geq1$}.

For \zw{ $0<t< 1$}, separate $A_{n}$ into the following two parts.
\begin{eqnarray}\label{more21}
\nonumber
A_{n1}&=& \int_{4 \log b_n}^{2(\frac{1}{t}-1)b_n^2} \exp\left(\frac{b_n^2}{2}\left( 1-(1+\pzx{tz} a_n^2)^{2/t}\right)\right) (1+\pzx{tz} a_n^2)^{1/t-1} dz
\\
 \nonumber
&<& \int_{4 \log b_n}^{2(\frac{1}{t}-1)b_n^2} e^{-z} (1+2(1-t))^{\frac{1}{t}-1}dz
\\
&=& O(b_n^{-4})
\end{eqnarray}
since
$\exp\left(\frac{b_n^2}{2}\left( 1-(1+\pzx{tz} a_n^2)^{2/t}\right)\right)<e^{-z}$.
For the second part,
\begin{eqnarray}\label{more22}
\nonumber
A_{n2}&=& \int_{2({1}/{t}-1)b_n^2}^{\infty} \exp\left(\frac{b_n^2}{2}\left( 1-(1+\pzx{tz} a_n^2)^{2/t}\right)\right) (1+\pzx{tz} a_n^2)^{1/t-1} dz
\\
\nonumber
&<&  \int_{2({1}/{t}-1)b_n^2}^{\infty} e^{-z}(t za_n^2)^{1/t-1}
\left(1+\frac{1}{\pzx{tz} a_n^2 } \right)^{1/t-1} dz
\\
\nonumber
&<& (ta_n^2 )^{1/t-1}  \left(1+\frac{1}{2(1-t)} \right)^{1/t-1}
\int_{2({1}/{t}-1)b_n^2}^{\infty} e^{-z}z^{1/t-1} dz
\\
\nonumber
&=& 2 (3-2t)^{1/t-1} e^{-2(1/t-1)b_n^2}
\\
&=& o(b_n^{-4}).
\end{eqnarray}
Hence,  \eqref{more21} and \eqref{more22} shows that \eqref{more2} holds as \pzx{$0<t<1$}.

Now changing to the case of \pzx{$t\ge1$}, by using Mills' inequality we have
\begin{eqnarray}\label{more23}
A_{n}&=& \int_{ 4 \log b_n}^{\infty} \exp\left( \frac{b_n^2}{2} \right)
\exp\left( -\frac{b_n^2 (1+\frac{tz}{b_n^2})^{2/t}}{2} \right) (1+\frac{tz}{b_n^2})^{1/t-1}dz
\nonumber\\
&=& b_n \exp\left( \frac{b_n^2}{2} \right)
\int_{b_n \left(1+\frac{4t\log b_n}{b_n^2}\right)^{1/t}}^{\infty} \exp\left(-\frac{s^2}{2}\right)ds
\nonumber\\
&=& \sqrt{2\pi} b_n \exp\left( \frac{b_n^2}{2} \right)
 \left(1-\Phi\left( b_n\left(1+\frac{4t\log b_n}{b_n^2} \right)^{1/t} \right)  \right)
\nonumber\\
&<& \frac{\exp\left( \frac{b_n^2}{2} \left( 1-\left(1+\frac{4t\log b_n}{b_n^2} \right)^{2/t}\right) \right)}
{\left(1+\frac{4t\log b_n}{b_n^2} \right)^{1/t}}
\nonumber\\
&<& \frac{\exp \left( -4\log b_n+\frac{8(t-2)(\log b_n)^2}{b_n^2} \right) }
{\left(1+\frac{4t\log b_n}{b_n^2} \right)^{1/t}}
\nonumber\\
&=& O(b_n^{-4})
\end{eqnarray}
since $(1+s)^{\frac{2}{t}} \geq 1+\frac{2}{t}s+\frac{1}{t}\left(\frac{2}{t}-1\right)s^2$ for $s>0$.

Combining with \eqref{more1}-\eqref{more2},
the proof of \eqref{add5} is complete. \qed

In order to show the second order asymptotic expansions of extreme value distributions, let
 \begin{eqnarray}\label{eq3.3aaa}
 \tilde{\Delta}(F_{\rho_n,t}^{n},H_\lambda;\pzx{c_{n}, d_{n}};x,y)=
 \P\left( M_{n1}\leq \omega_{n,t}(x),M_{n2}\leq \omega_{n,t}(y) \right)- H_{\lambda}(x,y).
 \end{eqnarray}

\begin{lemma}\label{lemma2}
Assume that the conditions of Theorem \ref{thm2} hold. Then,
 \begin{eqnarray}
\label{eq3.3}
\lim_{n\to\infty}(\log n)\tilde{\Delta}(F_{\rho_n,t}^{n},H_\lambda;\pzx{c_{n}, d_{n}};x,y)
=\frac{1}{2}\tau(\alpha,\lambda,x,y,t)H_\lambda(x,y)
\end{eqnarray}
where $\tau(\alpha,\lambda,x,y,t)$ is the one given in Theorem \ref{thm2}.
\end{lemma}

\noindent
\textbf{Proof.}~~ By using \eqref{eq3.2} and \eqref{add5}, we have
\begin{eqnarray*}
&&n\P\left(X>\omega_{n,t}(x),Y>\omega_{n,t}(y)\right)
\\
&=& \zw{\bar{\Phi}_{n,t}(y)}-\int_{y}^{\infty}\Phi\left(\frac{\omega_{n,t}(x)-\rho_n \omega_{n,t}(z)}{\sqrt{1-\rho_n^2}}\right)
e^{-z} \left( 1+\left((1-t)z-\frac{2-t}{2}z^2\right)b_n^{-2} \right)dz+O(b_n^{-4}).
\end{eqnarray*}
 for large $n$. It follows from \eqref{eq3.1b} and \eqref{add4} that
 \begin{eqnarray}\label{eq3.4}
&&
\nonumber
b_n^2 \int_{y}^{\infty}\left(\lambda+\frac{x-z}{2\lambda}-
\frac{\omega_{n,t}(x)-\rho_n \omega_{n,t}(z)}{\sqrt{1-\rho_n^2}}\right) \varphi\left(\lambda+\frac{x-z}{2\lambda}\right)e^{-z}dz
\\
\nonumber
&\to&
\left(\alpha-\zw{\frac{1}{2}\lambda^3-\frac{1}{2}\alpha\lambda^{-2}x-\frac{1}{4}\lambda x }
-\frac{1-t}{4\lambda}x^2\right)I_0
-\left(\frac{3}{4}\lambda-\frac{1}{2}\alpha\lambda^{-2}\right)I_1
+\frac{1-t}{4\lambda}I_2
\\
&=&
 \kappa_1(\alpha,\lambda,x,y,t)
\end{eqnarray}
as $n\to\infty$, where $I_{k}$ is the one given by \eqref{eq3.1a} and
\begin{eqnarray*}
&& \kappa_1(\alpha,\lambda,x,y,t)
\\
&=&
2\Big(\zw{(2-t)\lambda^4-(2-t)\lambda^2 x+(1-t)\lambda^2}\Big)\bar\Phi\left(\lambda+\frac{y-x}{2\lambda}\right)e^{-x}
\\
&& +\Big(2\alpha-\zw{(5-2t)\lambda^3+(1-t)\lambda x +(1-t)\lambda y}\Big)
\varphi\left(\lambda+\frac{y-x}{2\lambda}\right)e^{-x}.
\end{eqnarray*}

Note that by Taylor's expansion with Lagrange remainder term,
\begin{eqnarray}\label{add6}
&&
\Phi\left(\frac{\omega_{n,t}(x)-\rho_n \omega_{n,t}(z)}{\sqrt{1-\rho_n^2}}\right)\nonumber
\\
&=& \Phi\left(\lambda+\frac{x-z}{2\lambda}\right)+
\varphi\left(\lambda+\frac{x-z}{2\lambda}\right)
\left(\frac{\omega_{n,t}(x)-\rho_n \omega_{n,t}(z)}{\sqrt{1-\rho_n^2}}- \lambda-\frac{x-z}{2\lambda}\right)\nonumber
\\
&& +\frac{1}{2}v_n \varphi(v_n)
\left(\frac{\omega_{n,t}(x)-\rho_n \omega_{n,t}(z)}{\sqrt{1-\rho_n^2}}- \lambda-\frac{x-z}{2\lambda}\right)^2,
\end{eqnarray}
where $v_{n}$ is between $\frac{\omega_{n,t}(x)-\rho_n \omega_{n,t}(z)}{\sqrt{1-\rho_n^2}}$ and
 $\lambda+\frac{x-z}{2\lambda}$. By arguments similar to \eqref{eq3.4}, one can check that
\begin{eqnarray}\label{add7}
\int_{y}^{\infty} \left(\frac{\omega_{n,t}(x)-\rho_n \omega_{n,t}(z)}{\sqrt{1-\rho_n^2}}- \lambda-\frac{x-z}{2\lambda}\right)^2 v_n \varphi(v_n)e^{-z}dz
=O(b_n^{-4})
\end{eqnarray}
holds for large $n$. Hence from \eqref{eq3.4}, \eqref{add6} and \eqref{add7}, it follows that
\begin{eqnarray}
\label{eq3.5}
\lim_{n \to \infty}b_n^2 \int_{y}^{\infty}
\left(\Phi\left(\lambda+\frac{x-z}{2\lambda}\right)-
\Phi\left(\frac{\omega_{n,t}(x)-\rho_n \omega_{n,t}(z)}{\sqrt{1-\rho_n^2}}\right)\right)e^{-z}dz
=\kappa_1(\alpha,\lambda,x,y,t).
\end{eqnarray}

Note that
\begin{eqnarray}\label{eq3.10a}
\zw{\bar{\Phi}_{n,t}(x)}=e^{-x}-b_n^{-2}\mu(x)+O(b_n^{-4}),
\end{eqnarray}
cf. Theorem 1 in Hall (1980).
Now combining with \eqref{eq3.2}, \eqref{eq3.5} and \eqref{eq3.10a}, we have
\begin{eqnarray*}
&& b_n^2 \Big[\P\left( M_{n1}\leq \omega_{n,t}(x),M_{n2}\leq \omega_{n,t}(y) \right)- H_{\lambda}(x,y)\Big]\\
\\
&=&
b_n^2 H_{\lambda}(x,y) (1+o(1)) \Bigg[ - \bar{\Phi}_{n,t}(x)- \bar{\Phi}_{n,t}(y)+n\P(X>\omega_{n,t}(x),Y>\omega_{n,t}(y))
\\
&&
+\Phi\left(\lambda+\frac{x-y}{2\lambda}\right)e^{-y}+\Phi\left(\lambda+\frac{y-x}{2\lambda}\right)e^{-x} \Bigg]
\\
&=&
b_n^2 H_{\lambda}(x,y)(1+o(1)) \Bigg[ -\bar{\Phi}_{n,t}(x)+e^{-x} +
\int_{y}^{\infty} \Phi\left(\lambda+\frac{x-z}{2\lambda}\right)e^{-z}dz
\nonumber
\\
&&
-\int_{y}^{\infty} \Phi\left(\frac{\omega_{n,t}(x)-\rho_n \omega_{n,t}(z)}{\sqrt{1-\rho_n^2}}\right) e^{-z}
\left(1+\left((1-t)z-\frac{2-t}{2}z^2\right)b_n^{-2}\right)dz + O(b_n^{-4}) \Bigg]
\\
&\to&
H_{\lambda}(x,y) \Big[ \mu(x)+\kappa_1(\alpha,\lambda,x,y,t)-\kappa_2(\alpha,\lambda,x,y,t) \Big]
\\
&=&
H_{\lambda}(x,y) \tau(\alpha,\lambda,x,y,t)
\end{eqnarray*}
as $n\to\infty$, where
where
\begin{eqnarray*}
&&
\kappa_2(\alpha,\lambda,x,y,t)\\
&=&\int_{y}^{\infty} \Phi\left(\lambda+\frac{x-z}{2\lambda}\right)e^{-z}
\left((1-t)z-\frac{2-t}{2}z^2\right)dz
\\
&=& -\left(\frac{2-t}{2}y^2+y+1\right)\Phi\left(\lambda+\frac{x-y}{2\lambda}\right)e^{-y}
\\
&&
+ \left( 2 (2-t)\lambda^4-2(2-t)\lambda^2 x+2(1-t)\lambda^2+\frac{2-t}{2}x^2+x+1 \right)
\bar{\Phi}\left(\lambda+\frac{y-x}{2\lambda}\right)e^{-x}
\\
&&
\pzx{-  \Big( 2(2-t)\lambda^3-(2-t)\lambda (x+y)-2\lambda \Big)}
\varphi\left(\lambda+\frac{y-x}{2\lambda}\right)e^{-x}
\end{eqnarray*}
and $\tau(\alpha,\lambda,x,y,t)$ is the one given by Theorem \ref{thm2}.
The proof is complete. \qed

\begin{lemma}  \label{lemma3}
With $c_n$ and $d_n$ given by \eqref{eq1.7}, the following results hold.
 \begin{itemize}
\item[(a).]
For the case of $\lambda=\infty$,

(i) if $\rho_n \in [-1,0]$, we have
\begin{eqnarray}\label{eq3.8}
\lim_{n\to \infty}(\log n)\tilde\Delta(F_{\rho_n,t}^{n},H_{\infty};\pzx{c_{n}, d_{n}};x,y)
=\frac{1}{2}\Big(\mu(x)+\mu(y)\Big)  H_{\infty}(x,y).
\end{eqnarray}
(ii) if $\rho_n \in (0,1)$ and $\lim_{n \to \infty} \frac{\log b_n}{b_n^2(1-\rho_n)}=0$,
then \eqref{eq3.8} also holds.

\item[(b).]
For the case of $\lambda=0$,

(i) if $\rho_n=1$, we have
\begin{eqnarray}\label{eq3.9}
\lim_{n\to \infty}(\log n)\tilde\Delta(F_{\rho_n,t}^{n},H_0;\pzx{c_{n}, d_{n}};x,y)
=\frac{1}{2}\mu \Big(\min(x,y)\Big)  H_0(x,y).
\end{eqnarray}
(ii) if $\rho_n \in (0,1)$ and $\lim_{n \to \infty}b_n^{6}(1-\rho_n)=c_1 \in [0,\infty)$, then
\eqref{eq3.9} also holds.
\end{itemize}
\end{lemma}

\noindent
\textbf{Proof.}~~ For $\lambda=\infty$, we first consider (i), i.e., the case of $\rho_n \in [-1,0]$. Note that either complete independent $(\rho_n\equiv0)$ or complete negative dependent $(\rho_n\equiv-1)$ both imply
$\lambda=\infty$. Thus from \eqref{eq3.10a} it follows that both
\begin{eqnarray}
\label{add91}
\nonumber
&&  b_n^2 \left( -n(1-F_{-1}(\omega_{n,t}(x),\omega_{n,t}(y)))+e^{-y}+e^{-x} \right)
\\
\nonumber
&=&
 b_n^2\left( - \bar{\Phi}_{n,t}(x)+e^{-x} \right) +
b_n^2\left( - \bar{\Phi}_{n,t}(y)+e^{-y} \right)
  + nb_n^2\P (\omega_{n,t}(x)<X<-\omega_{n,t}(y))
  \\
&\to &  \mu(x)+\mu(y)
\end{eqnarray}
and
\begin{eqnarray}
\label{add92}
\nonumber
&&  b_n^2 \left( -n(1-F_{0}(\omega_{n,t}(x),\omega_{n,t}(y)))+e^{-y}+e^{-x} \right)
\\
\nonumber
&=&
 b_n^2\left( - \bar{\Phi}_{n,t}(x)+e^{-x} \right)
 +b_n^2 \left( - \bar{\Phi}_{n,t}(y)+e^{-y} \right)
+\frac{b_n^2}{n}\bar{\Phi}_{n,t}(x)\bar{\Phi}_{n,t}(y)
\\
&\to&  \mu(x)+\mu(y)
\end{eqnarray}
hold as $n\to\infty$,  showing that the claimed results \eqref{eq3.8} hold for  $\rho_n\equiv-1$ and  $\rho_n\equiv0$ respectively.
Thus, it follows from Slepian's Lemma that \eqref{eq3.8} also holds for $\rho_n \in [-1,0]$.

Now switch to the case of $\rho_n \in (0,1)$ with additional condition $\lim_{n \to \infty}\frac{\log b_n}{b_n^2(1-\rho_n)}=0$, implying $\lambda=\infty$.
For fixed $x,z \in \mathbb{R}$, one can check that
\begin{eqnarray}\label{add9}
\lim_{n \to \infty}\frac{\omega_{n,t}(x)-\rho_n \omega_{n,t}(z)}{\sqrt{1-\rho_n^2}}=\infty.
\end{eqnarray}
Note that the condition $\lim_{n \to \infty}\frac{\log b_n}{b_n^2(1-\rho_n)}=0$ implies
$\lim_{n \to \infty}b_n^2(1-\rho_n)=\infty$. By \eqref{add9} and Mills' inequality,
\begin{eqnarray}
\label{eq3.18}
&& b_n^4 \left( 1-\Phi\left(\frac{\omega_{n,t}(x)-\rho_n \omega_{n,t}(z)}{\sqrt{1-\rho_n^2}}\right) \right)
\nonumber
\\
&<& \frac{b_n^4  \exp\left(-\frac{(\omega_{n,t}(x)-\rho_n \omega_{n,t}(z))^2}{2(1-\rho_n^2)}\right)}
{ \frac{\omega_{n,t}(x)-\rho_n \omega_{n,t}(z)}{\sqrt{1-\rho_n^2}} }
\nonumber
\\
&=&  \frac{\exp\left(-\frac{((c_n x+d_n)^{1/t}-\rho_n (c_n z+d_n)^{1/t})^2}{2(1-\rho_n^2)}+4\log b_n\right)}
{ \frac{(c_n x+d_n)^{1/t}-\rho_n (c_n z+d_n)^{1/t}}{\sqrt{1-\rho_n^2}}}
\nonumber
\\
&=&
 \left( b_n\sqrt{\frac{1-\rho_n}{1+\rho_n}}+\frac{x-z}{b_n\sqrt{1-\rho_n^2}}+
 \frac{z}{b_n}\sqrt{\frac{1-\rho_n}{1+\rho_n}}
+\frac{(1-t)(x^2-\rho_n z^2)}{2b_n^3\sqrt{1-\rho_n^2}} +\sqrt{\frac{1-\rho_n}{1+\rho_n}}O(b_n^{-5})\right)^{-1}
\nonumber
\\
&&  \times  \exp\left( -\frac{b_n^2(1-\rho_n)}{2(1+\rho_n)}-\frac{(x-\rho_n z)^2}{2b_n^2(1-\rho_n^2)}-
\frac{(1-t)^2 (x^2-\rho_n z^2)^2}{8b_n^6(1-\rho_n^2)}-\frac{x-\rho_n z}{1+\rho_n}-\frac{(1-t)(x^2-\rho_n z^2)}{2b_n^2(1+\rho_n)} \right.
\nonumber
\\
&& \left.
-\frac{(x-\rho_n z)(1-t)(x^2-\rho_n z^2)}{2b_n^4(1-\rho_n^2)} +\frac{1-\rho_n}{2(1+\rho_n)}O(b_n^{-5})+4\log b_n \right)
\nonumber
\\
&<& (1+o(1))  e^{-\frac{x-z}{2}}
\exp \left\{ -\frac{b_n^2(1-\rho_n)}{2(1+\rho_n)} \left( 1- \frac{8(1+\rho_n)\log b_n}{b_n^2(1-\rho_n)}
+\frac{(1+\rho_n)\log b_n^2(1-\rho_n)}{b_n^2 (1-\rho_n)}\right) \right\}
\nonumber
\\
& \to & 0
\end{eqnarray}
as $n \to \infty$.
Note that
\begin{eqnarray}\label{eq3.11a}
n^{-1}=\pzx{\bar{\Phi}(b_n)}=\frac{\varphi(b_n)}{b_n}(1-b_n^{-2}+O(b_n^{-4})).
\end{eqnarray}
Hence, by using \eqref{more1}-\eqref{more2}, we have
\begin{eqnarray}
\label{eq3.11b}
&& n \P(X>\omega_{n,t}(x),Y>\omega_{n,t}(y))
\nonumber
\\
\nonumber
&=& b_n^{-4}(1-b_n^{-2}+O(b_n^{-4}))^{-1}
\int_{y}^{\infty} b_n^4 \bar{\Phi}(\frac{\omega_{n,t}(x)-\rho_n \omega_{n,t}(z)}{\sqrt{1-\rho_n^2}})
\exp\left(\frac{b_n^2}{2}\left( 1-(1+tza_n^2)^{2/t}\right)\right)
\\
\nonumber
&& \times (1+tza_n^2)^{1/t-1}dz
\\
&=& O(b_n^{-4}).
\end{eqnarray}
 It follows from \eqref{eq3.10a} and \eqref{eq3.11b} that
 \begin{eqnarray}
 \label{eq3.12a}
&& b_n^2 \Big[F_{\rho_n}^n(\omega_{n,t}(x),\omega_{n,t}(y))-H_{\infty}(x,y)\Big]
\nonumber
\\
\nonumber
&=&
b_n^2 H_{\infty}(x,y) (1+o(1)) \Big[-n (1-F_{\rho_n}(\omega_{n,t}(x),\omega_{n,t}(y)))+e^{-x}+e^{-y}\Big]
\\
\nonumber
&=&
b_n^2 H_{\infty}(x,y) (1+o(1)) \Big[- (e^{-x}+e^{-y}-b_n^{-2}(\mu(x)+\mu(y)+O(b_n^{-2})))+e^{-x}+e^{-y}  \Big]
\\
\nonumber
&\to&  H_{\infty}(x,y) \Big[\mu(x)+\mu(y)\Big]
\end{eqnarray}
as $n\to\infty$. Proof the case of $\lambda=\infty$ is complete.

(b).
For the case of $\lambda=0$, we first consider the complete positive dependence case $(\rho_n\equiv1)$.  Without loss of generality, assume that
$y<x$, we have
\begin{eqnarray}
\label{eq3.21a}
b_n^2 \Big[-n(1-F_1(\omega_{n,t}(x),\omega_{n,t}(y)))+e^{-y}\Big]
= b_n^2\Big[-\zw{\bar{\Phi}_{n,t}(y)}+e^{-y}\Big]
\to\mu(y)
\end{eqnarray}
as $n \to \infty$ since \eqref{eq3.10a} holds.
The rest is for the case of $\rho_n \in (0,1)$.
For $y<x \in \R$, if
$\max (x,y)=x<z<4 \log b_n$ we have
 \begin{eqnarray}\label{add10}
&& \Phi\left(\frac{\omega_{n,t}(y)-\rho_n \omega_{n,t}(z)}{\sqrt{1-\rho_n^2}}\right)
\nonumber\\
&<& -\frac{\varphi\left(\frac{\omega_{n,t}(y)-\rho_n \omega_{n,t}(z)}{\sqrt{1-\rho_n^2}}\right)}
{\frac{\omega_{n,t}(y)-\rho_n \omega_{n,t}(z)}{\sqrt{1-\rho_n^2}}}
\nonumber\\
&=&
\frac{\exp \left( -\frac{1}{2} \left(\lambda_n+\frac{y-z}{2\lambda_n}\right)^2 (1+o(1)) \right)}
{ \left( - \lambda_n+\frac{z-y}{2\lambda_n} \left( 1+\frac{(1-t)(y+z)}{2 b_n^2} \right) -\frac{\lambda_n z}{b_n^2} -
\frac{(1-t)\lambda_n z^2}{2b_n^{4}} +
\lambda_n O(b_n^{-6})  \right) (1-\frac{\lambda_n^2}{b_n^2})^{-\frac{1}{2}}  }
\nonumber\\
&=&
\frac{\exp \left(-\frac{1}{2} \left(\lambda_n+\frac{y-z}{2\lambda_n}\right)^2 (1+o(1)) \right)}
{\frac{z-y}{2\lambda_n}(1+o(1))}
\end{eqnarray}
for large $n$ due to $\Phi(-x)=\pzx{\bar{\Phi}(x)}$ and Mills' inequality since
$\frac{\omega_{n,t}(y)-\rho_n \omega_{n,t}(z)}{\sqrt{1-\rho_n^2}}<0$ for large $n$
when $\lim_{n \to \infty} b_n^{6} (1-\rho_n)=c_1 \in [0,\infty)$. Therefore,
\begin{eqnarray}
\label{eq3.12}
\nonumber
&& \int_{x}^{4\log b_n}\Phi\left( \frac{\omega_{n,t}(y)-\rho_n \omega_{n,t}(z)}{\sqrt{1-\rho_n^2}} \right)
\exp\left(\frac{b_n^2}{2}\left( 1-(1+\pzx{tz} a_n^2)^{2/t}\right)\right) (1+\pzx{tz} a_n^2)^{1/t-1}dz
\\
\nonumber
&<&
\frac{2\lambda_n}{x-y}(1+o(1))  \int_{x}^{4\log b_n}
\exp \left( -\frac{\lambda_n^2}{2}-\frac{y-z}{2}-\frac{(y-z)^2}{8\lambda_n^2}-z+o(b_n^{-1})+
\left( \frac{1}{t}-1\right)\log (1+\pzx{tz} a_n^2) \right)dz
\\  \nonumber
&=&
2\lambda_n (1+o(1)) \frac{ \exp \left( -\frac{\lambda_n^2}{2}-\frac{y}{2}-\frac{(y-x)^2}{8\lambda_n^2} \right)}  {x-y}
\int_{x}^{4\log b_n} \exp \left(-\frac{z}{2}\right)dz
\\  \nonumber
&<&
4\lambda_n b_n^{-2}(1+o(1)) \frac{\exp \left( -\frac{\lambda_n^2}{2}-\frac{y}{2}-\frac{(y-x)^2}{8\lambda_n^2} \right)}{y-x}
\\
&=& O(b_n^{-4})
\end{eqnarray}
for large $n$ by using $\lim_{n \to \infty}b_n^{6}(1-\rho_n)=c_{1}$. \pzx{It follows from \eqref{more2} that}
\begin{eqnarray}
\label{eq3.13}
\nonumber
&& \int_{4\log b_n}^{\infty} \Phi\left( \frac{(\omega_{n,t}(y)-\rho_n \omega_{n,t}(z))^2}{\sqrt{1-\rho_n^2}} \right)
\exp\left(\frac{b_n^2}{2}\left( 1-(1+tza_n^2)^{2/t}\right)\right) (1+tza_n^2)^{1/t-1}dz
\\  \nonumber
&<& \int_{4\log b_n}^{\infty} \exp\left(\frac{b_n^2}{2}\left( 1-(1+tza_n^2)^{2/t}\right)\right)
(1+tza_n^2)^{1/t-1}dz
\\
&=& O(b_n^{-4}).
\end{eqnarray}
Combining \eqref{eq3.10a}, \eqref{eq3.12} and \eqref{eq3.13}, for $y<x$ we have
\begin{eqnarray*}
&& 1-F_{\rho_n}\Big(\omega_{n,t}(\min(x,y)),\omega_{n,t}(\max(x,y))\Big)
\\
&=& \zw{\bar{\Phi}_{n,t}(y)+\bar{\Phi}_{n,t}(x)}
-\P\Big(X>\omega_{n,t}(y),Y>\omega_{n,t}(x)\Big)
\\
&=& \zw{\bar{\Phi}_{n,t}(y)+\bar{\Phi}_{n,t}(x)}
\\
&& -\int_{x}^{\infty}\bar{\Phi}\left(\frac{\omega_{n,t}(y)-\rho_n \omega_{n,t}(z)}{\sqrt{1-\rho_n^2}}\right)
\exp\left(\frac{b_n^2}{2}\left( 1-(1+tza_n^2)^{\frac{2}{t}}\right)\right) (1+tza_n^2)^{\frac{1}{t}-1}dz
\\
&=& \zw{\bar{\Phi}_{n,t}(y)}+n^{-1}(1-b_n^{-2}+O(b_n^{-4}))^{-1}
\\
&& \times \int_{x}^{\infty} \Phi\left(\frac{\omega_{n,t}(y)-\rho_n \omega_{n,t}(z)}{\sqrt{1-\rho_n^2}}\right)
\exp\left(\frac{b_n^2}{2}\left( 1-(1+tza_n^2)^{\frac{2}{t}}\right)\right) (1+tza_n^2)^{\frac{1}{t}-1}dz
\\
&=& n^{-1}(e^{-y}-b_n^{-2}\mu(y)+O(b_n^{-4}))
\end{eqnarray*}
holds for large $n$, which implies the desired result. The proof is complete.
\qed

\begin{lemma}  \label{lemma4}
For $t=2$, with $c_{n}^{*}$ and $d_{n}^{*}$ given by \eqref{add2}, the following results hold.
 \begin{itemize}
\item[(a).]
For the case of $\lambda=\infty$,

(i) if $\rho_n \in [-1,0]$, we have
\begin{eqnarray}\label{eq3.14}
\lim_{n\to \infty}(\log n)^2\tilde\Delta(F_{\rho_n,2}^{n},H_{\infty};\pzx{c_{n}^{*}, d_{n}^{*}};x,y)
=\frac{1}{4}\Big(\nu(x)+\nu(y)\Big)  H_{\infty}(x,y).
\end{eqnarray}
(ii) if $\rho_n \in (0,1)$ and $\lim_{n \to \infty} \frac{\log b_n}{b_n^2(1-\rho_n)}=0$,
then \eqref{eq3.14} also holds.

\item[(b).]
For the case of $\lambda=0$,

(i) if $\rho_n=1$, we have
\begin{eqnarray}\label{eq3.15}
\lim_{n\to \infty}(\log n)^2\tilde\Delta(F_{\rho_n,2}^{n},H_0;\pzx{c_{n}^{*}, d_{n}^{*}};x,y)
=\frac{1}{4}\nu\Big(\min(x,y)\Big)  H_0(x,y).
\end{eqnarray}
(ii) if $\rho_n \in (0,1)$ and $\lim_{n \to \infty}b_n^{14}(1-\rho_n)=c_2 \in [0,\infty)$, then \eqref{eq3.15} also holds.
\end{itemize}
\end{lemma}

\noindent
\textbf{Proof.}~~
(a). For $\lambda=\infty$.
Firstly note that
\begin{eqnarray}\label{eq3.16}
\zw{\bar{\Phi}_{n,2}(x)}=e^{-x}-b_n^{-4}\nu(x)+O(b_n^{-6})
\end{eqnarray}
derived by Theorem 1 in Hall (1980).
\zw{Arguments similar to that of \eqref{add91} and \eqref{add92}, by using \eqref{eq3.16} we have
\begin{eqnarray*}
b_n^4 \left( -n (1-F_{-1}(\omega_{n,2}^{*}(x),\omega_{n,2}^{*}(y))) +\Phi(\lambda+\frac{x-y}{2\lambda})e^{-y} +\Phi(\lambda+\frac{y-x}{2\lambda})e^{-x} \right)
\to  \nu(x)+\nu(y)
\end{eqnarray*}
as $n \to \infty$ for $\rho_n\equiv-1$, and
\begin{eqnarray*}
b_n^4 \left( -n (1-F_{0}(\omega_{n,2}^{*}(x),\omega_{n,2}^{*}(y))) +\Phi(\lambda+\frac{x-y}{2\lambda})e^{-y} +\Phi(\lambda+\frac{y-x}{2\lambda})e^{-x} \right)
\to  \nu(x)+\nu(y)
\end{eqnarray*}
also holds as $n \to \infty$ for $\rho_n\equiv0$}. Therefore, \eqref{eq3.14} holds for $\rho_n\equiv-1$ and $\rho_n\equiv0$ respectively.
By using Slepian's Lemma, \eqref{eq3.14} also holds for $\rho_n \in [-1,0]$.

Now switch to the case of $\rho_n \in(0,1)$ with additional condition $\lim_{n \to \infty} \frac{\log b_n}{b_n^2(1-\rho_n)}=0$, implying $\lambda=\infty$.
For fixed $x$ and $z$, note that
\begin{eqnarray*}
\frac{\omega_{n,2}^{*}(x)-\rho_n \omega_{n,2}^{*}(z)}{\sqrt{1-\rho_n^2}}
 \to \infty
\end{eqnarray*}
as $n\to\infty$ since $\lambda_n^2=\frac{b_n^2}{2}(1-\rho_n) \to \infty$ as $n \to \infty$. Therefore for $t=2$, arguments similar to \eqref{eq3.18}, we have
\begin{eqnarray*}
 b_n^6 \left( 1-\Phi\left(\frac{\omega_{n,2}^{*}(x)-\rho_n \omega_{n,2}^{*}(z)}{\sqrt{1-\rho_n^2}}\right) \right)
\to  0
\end{eqnarray*}
as $n \to \infty$. Hence, it follows from \eqref{eq3.11a} that
\begin{eqnarray*}
&& \P(X>\omega_{n,2}^{*}(x),Y>\omega_{n,2}^{*}(y))
\\
&=& n^{-1} b_n^{-6}(1+b_n^{-2}+O(b_n^{-4}))
\int_{y}^{\infty} b_n^6
 \left(1-{\Phi}(\frac{\omega_{n,2}(x)-\rho_n \omega_{n,2}(y)}{\sqrt{1-\rho_n^2}})\right)
 \\
&& \times
\exp\left(-z+(1+z)a_n^2\right) \left(1-a_n^2\right) \left(1+2\left(z-(1+z)a_n^2\right)a_n^2\right)^{-\frac{1}{2}} dz
\\
&=& O(n^{-1}b_n^{-6})
\end{eqnarray*}
for large $n$. Hence,
 \begin{eqnarray*}
 b_n^4 \Big[F_{\rho_n}^n(\omega_{n,t}^{*}(x),\omega_{n,t}^{*}(y))-H_{\infty}(x,y)\Big]
\to  H_{\infty}(x,y) \Big[\nu(x)+\nu(y)\Big]
\end{eqnarray*}
holds as $n\to\infty$.
The proof the case of $\lambda=\infty$ is complete.

(b).
For the case of $\lambda=0$. We first consider the complete positive dependence case $(\rho_n\equiv1)$, without loss
of generality, assume that $y<x$. Hence,
\begin{eqnarray*}
 b_n^{4}\Big[ -n(1-F_{1}(\omega_{n,2}^{*}(x),\omega_{n,2}^{*}(y)))+\Phi(\lambda+\frac{x-y}{2\lambda})e^{-y}+\Phi(\lambda+\frac{y-x}{2\lambda})e^{-x} \Big]
\to  \nu(y)
\end{eqnarray*}
as $n\to\infty$ provided that \eqref{eq3.16} holds. The remainder is to prove the case of $\rho_n \in(0,1)$ with $\lim_{n \to \infty}b_n^{14}(1-\rho_n)=c_2 \in [0,\infty)$.
By arguments similar to  that of \eqref{eq2.12} and \eqref{eq2.13}, for fixed $y<x \in \R$ we have
\begin{equation}
\label{eq3.17}
\int_{x}^{\infty} \Phi \left(\frac{\omega_{n,2}^{*}(y)-\rho_n \omega_{n,2}^{*}(z)}{\sqrt{1-\rho_n^2}}\right)
\exp(-z+ (1+z)a_n^2) \frac{1-a_n^2}{(1+2 (z-(1+z)a_n^2)a_n^2)^{\frac{1}{2}}}  dz
= O(b_n^{-6})
\end{equation}
for large $n$ by using $\lim_{n \to \infty}b_n^{14}(1-\rho_n)=c_{2}$ . Combining \eqref{eq3.16} with \eqref{eq3.17}, we have
\begin{eqnarray*}
&& 1-F_{\rho_n}(\omega_{n,2}^{*}(x),\omega_{n,2}^{*}(y))
\\
&=& \zw{\bar{\Phi}_{n,2}(y)}+n^{-1}(1-b_n^{-2}+O(b_n^{-4}))^{-1}
\\
&& \times
\int_{x}^{\infty} \Phi \left(\frac{\omega_{n,2}^{*}(y)-\rho_n \omega_{n,2}^{*}(z)}{\sqrt{1-\rho_n^2}} \right)
\exp(-z+ (1+z)a_n^2) \frac{1-a_n^2}{(1+2 (z-(1+z)a_n^2)a_n^2)^{\frac{1}{2}}}   dz
\\
&=& n^{-1} \left( e^{-y}-b_n^{-4}\nu(y)+O(b_n^{-6})  \right)
\end{eqnarray*}
for large $n$, which implies the desired result. The proof is complete.
\qed

\section{Proofs}
\label{sec4}

\noindent \textbf{Proof of Theorem \ref{thm1}.}~~
Obviously,
 \begin{eqnarray*}
 &&\P\left( \abs{M_{n1}}^t\leq c_n x+d_n, \abs{M_{n2}}^t\leq c_n y+d_n \right)
 \\
 &=& F_{\rho_n,t}^{n}(\omega_{n,t}(x),\omega_{n,t}(y))-F_{\rho_n,t}^{n}(\omega_{n,t}(x),-\omega_{n,t}(y))
 -F_{\rho_n,t}^{n}(-\omega_{n,t}(x),\omega_{n,t}(y))+F_{\rho_n,t}^{n}(-\omega_{n,t}(x),-\omega_{n,t}(y)).
 \end{eqnarray*}
Note that
 \begin{eqnarray}\label{eq4.1}
 \nonumber
&& F_{\rho_n,t}^{n}(\omega_{n,t}(x),-\omega_{n,t}(y))+F_{\rho_n,t}^{n}(-\omega_{n,t}(x),\omega_{n,t}(y))
-F_{\rho_n,t}^{n}(-\omega_{n,t}(x),-\omega_{n,t}(y))
\\
\nonumber
&\leq& \P\left(M_{n2}\leq -\omega_{n,t}(y)\right)+\P\left(M_{n1}\leq -\omega_{n,t}(x)\right)
-\min \{  \Phi^n(-\omega_{n,t}(x)),\Phi^n(-\omega_{n,t}(y)) \}
\\
\nonumber
&=& \Phi^n(-\omega_{n,t}(x))+\Phi^n(-\omega_{n,t}(y))-\min \{ \Phi^n(-\omega_{n,t}(x)),\Phi^n(-\omega_{n,t}(y)) \}
\\
&=& o(b_n^{-4})
\end{eqnarray}
since
\begin{eqnarray*}
\bar{\Phi}^{n-1}(-\omega_{n,t}(x))=\left(n^{-1}e^{-x}(1+O(b_n^{-2}))\right)^{n-1}
=o(b_n^{-4}),
\end{eqnarray*}
cf. Lemma 3.1 in Zhou and Ling (2016).
Combining \eqref{eq4.1} with Lemma \ref{lemma1}, we can get the claimed result \eqref{eq2.2}. \qed

\noindent \textbf{Proof of Theorem \ref{thm2}.}~~
It follows from \eqref{eq4.1}  and Lemma \ref{lemma2} that
\begin{eqnarray*}
 \Delta(F_{\rho_n,t}^{n},H_{\lambda};\pzx{c_{n},d_{n}}; x,y)
=  \tilde{\Delta}(F_{\rho_n,t}^{n},H_{\lambda};\pzx{c_{n},d_{n}}; x,y)+ o(b_n^{-4}),
\end{eqnarray*}
 so the result \eqref{eq2.3} is obtained. \qed

\noindent \textbf{Proof of Theorem \ref{thm3} and Theorem \ref{thm4}.}~~
By using \eqref{eq4.1},  Lemma \ref{lemma3}  and Lemma \ref{lemma4} respectively to derive the desired results. \qed

\vspace{3cm}

\noindent {\bf Acknowledgments}~~ The first author was supported by the Fundamental Research Funds for the Central
Universities Grant no. XDJK2016E117, the CQ Innovation Project for Graduates Grant no. CYS16046.

\end{document}